\documentclass[12pt]{amsart}
\usepackage{txfonts}
\usepackage{bbm}
\usepackage{mathrsfs}
\usepackage{amssymb}
\usepackage{amssymb,amsmath,amsthm,amscd}
\usepackage[all]{xy}

\numberwithin{equation}{section}

\newtheoremstyle{fancy1}{10pt}{10pt}{\itshape}{12pt}{\textsc\bgroup}{.\egroup}{8pt}{
}
\newtheoremstyle{fancy2}{10pt}{10pt}{}{12pt}{\itshape}{.}{8pt}{ }

\theoremstyle{fancy1}

\newtheorem{cor}[equation]{Corollary}
\newtheorem{lem}[equation]{Lemma}

\newtheorem{thm}[equation]{Theorem}
\newtheorem{problem}[equation]{Problem}
\newtheorem{main}{Theorem}
\newtheorem*{main*}{Theorem}

\newtheorem*{cor*}{Corollary}

\theoremstyle{fancy2}

\newtheorem{rem}[equation]{Remark}
\newtheorem*{rem*}{Remark}
\newtheorem{example}[equation]{Example}

\newcommand{\cref}[1]{Corollary~\ref{#1}}

\newcommand{\Sph}{\mathbb{S}}
\newcommand{\Disc}{\mathbb{D}}

\newcommand{\R}{{\mathbb{R}}}

\newcommand{\G}{\ensuremath{\operatorname{\mathsf{G}}}}

\def\con#1=#2(#3){#1 \equiv #2 \bmod{#3}}

\newcommand{\diam}{\ensuremath{\operatorname{diam}}}
\renewcommand{\sec}{\ensuremath{\operatorname{sec}}}

\newcommand{\ric}{\ensuremath{\operatorname{ric}}}
\newcommand{\vol}{\ensuremath{\operatorname{vol}}}

\newcommand{\inj}{\ensuremath{\operatorname{inj}}}

\newcommand{\cd}{\ensuremath{\operatorname{cd}}}
\newcommand{\excess}{\ensuremath{\operatorname{excess}}}

\newcommand{\no}{\noindent}

\begin{document}
\date{\today}

\title{Manifolds tightly covered by two metric balls}


\author{Jianming Wan}
\address{School of Mathematics\\Northwest University\\
Xi'an 710127\\China } \email{wanj\_m@aliyun.com}

\thanks{The author was supported in part by National Natural Science Foundation of China
 No.11601421 and No.11601422}

\begin{abstract}
In this note we provide natural optimal geometric conditions for a Riemannian manifold suitably covered by two open metric balls to be homeomorphic to a sphere. This can be viewed as a geometric analogue of Brown's theorem in topology stating that a closed manifold covered by two topological balls is a sphere.
\end{abstract}

\maketitle



As the simplest closed manifold, the sphere enjoys a unique and basic role topologically as well as metrically. Geometrically, the unit sphere is uniquely determined as an ''optimal object" in a variety of ways often referred to as \emph{Sphere Theorems}. Examples of such recognition results include the classical Rauch-Berger-Klingenberg $1/4$-pinching theorem (for diffeomorphism see Brendle-Schoen \cite{[BS]}), the diameter sphere theorem \cite{[GS]}, Micallef-Moore's positive isotropic curvature sphere theorem \cite{[MM]}, and Perelman's almost maximal volume sphere theorem \cite{[P1]} (for diffeomorphism see Colding and Cheeger \cite{Co,CC}). Topologically, Brown's Theorem \cite{[B]} recognizes the sphere as the only closed manifold covered by two open euclidean balls.

As a metric contrast to Brown's Theorem, we point out, that any closed (smooth) manifold, $M$ admits a Riemannian metric so that it is covered by two (proper) open metric balls, even \emph{tightly covered} in the following sense:

For any $\epsilon >0$ and fixed $r>0$, there is a Riemannan metric on $M$ so that
$$M = B(p,r+\epsilon) \cup B(q,r+\epsilon), \ \text{with} \  \rho(p,q) = 2r$$
where $B(p,r)$ denotes the open $r$-ball centered at $p$ and $\rho(p,q)$ is the distance between $p$ and $q$. For example, it can be arranged that the complement of an arbitrarily small metric ball in $M$ is a disc with constant curvature $1$.

However, if for a fixed Riemannian metric $M$ is $\epsilon$-tightly covered for every $\epsilon >0$, then of course
$$M = D(p,r) \cup D(q,r), \ \text{with} \  \rho(p,q) = 2r,$$
where $D(p,r)$ denotes the closed $r$ ball with center $p$. In this case all geodesics emanating from $p$ of length $2r$ are minimal and terminate at $q$ (see Lemma \ref{suspension}). In particular, such an $M$ is a topological sphere.

Our goal is to seek natural geometric conditions under which a Riemannian manifold tightly covered by two open metric balls, in the sense above, is a topological sphere. Our results hinge on the observation that for certain classes of metric spaces being tightly covered by two proper open metric balls is equivalent to having small excess in the sense of \cite{[GP]}. Here $\text{excess} M < \delta$ means there is a pair of points $p,q \in M$ such that for any $x \in M$, $$\rho(p,x) + \rho(x,q) - \rho(p,q) < \delta.$$ In this case, clearly $M$ is $\epsilon = \frac{\delta}{2}$ tightly covered in the above sense.

Indeed, we have (see section 1)

\begin{main}
Let $\mathcal{M}$ be a Gromov-Hausdorff precompact class of closed Riemannian manifolds for which any $X \in \bar{\mathcal{M}}$ is a non-branching geodesic metric space. Then for any $\epsilon >0$ there is a $\delta > 0$ such that $\text{excess}M < \epsilon$ if $M$ is $\delta$-tightly covered by two open balls, and vice versa.
\end{main}

Recall, that from the Bishop-Gromov relative volume comparison theorem it follows that the class of all closed $n$-manifolds $M$ with  Ricci curvature, $\ric M \ge (n-1)k$ and diameter $\diam M \le D$ is Gromov-Hausdorff precompact. The subclasses where the sectional curvature $\sec M \ge k$, or the injectivity radius
$\inj M \ge i$ are examples of $\mathcal{M}$ as above. In the first case because any limit object is an Alexandrov space, and in the second case the non-branching property was proved by Taylor in \cite{Ta}.

Appealing to the main theorems of \cite{[Pet]} and \cite{[GP]} this yields the following immediate corollaries

\begin{main}
For any real $k, r > 0, i> 0$ and integer $n \ge 2$ there is an $\epsilon_0 = \epsilon_0(k,r,i,n)$ such the following holds: Any closed Riemannian $n$-manifold $M$ with $\ric M \ge (n-1)k$, $\inj M \ge i$ and
$$M^{n}=B_{p}(r+\epsilon)\cup B_{q}(r+\epsilon), \quad \rho(p,q) = 2r$$
is homeomorphic to $\Sph^n$ if $\epsilon < \epsilon_0$.
\end{main}

If the condition $\inj M \ge i$ is relaxed to $\vol M \ge v$, the conclusion fails as, e.g., the examples due to Anderson \cite{An} shows. However, if at the same time $\ric M \ge (n-1)k$ is strengthened to $\sec M \ge k$ we have:

\begin{main}
For any real $k, r > 0, v> 0$ and integer $n \ge 2$ there is an $\epsilon_1 = \epsilon_1(k,r,v,n)$ such the following holds: Any closed Riemannian $n$-manifold $M$ with $\sec M \ge k$, $\vol M \ge v$ and
$$M^{n}=B_{p}(r+\epsilon)\cup B_{q}(r+\epsilon), \quad \rho(p,q) = 2r$$
is homeomorphic to $\Sph^n$ if $\epsilon < \epsilon_1$.
\end{main}

In these statements we have no explicit estimate for $\epsilon_i$. Likewise, we do not prove that the open metric balls $B(p,r+\epsilon)$ and $B(q,r+\epsilon)$
 in $M$ are homeomorphic to the euclidean n-ball. Although, Theorem A implies that Theorems B and C are equivalent to the main results in \cite{[Pet]} and \cite{[GP]} we present alternate short proofs.

 \medskip

 In contrast, if $2r = d=\diam M$ and $\vol M \ge v$ is strengthened to $\inj M \ge i$, we have a constructive proof that being $\epsilon$-tightly covered implies small excess. Thus by critical point theory lemma 3 of \cite{[GP]} we have

\begin{main}
There is an explicit $\epsilon_2 = \epsilon_2(k,i,d)$, such that any closed Riemannian $n$-manifold $M$ with $\sec M \ge k$, $\inj M \ge i$ and
$$M^{n}=B_{p}(r+\epsilon)\cup B_{q}(r+\epsilon), \quad \rho(p,q) = 2r = d=\diam M$$
 is homeomorphic to $n$-sphere $\Sph^{n}$ if $\epsilon < \epsilon_2$.
\end{main}

Note also that here $\epsilon_2$ does not depend on $n$ either. When, e.g., $k=0$ it can be shown using critival point theorey in conjunction with Toponogov's comparison theorem that we can choose $\epsilon_2(0,i,d) = \frac{1}{2}(-d+\sqrt{d^{2}+i^{2}/2})$
\smallskip

 For basic tools and results in Riemannian and Alexandrov Geometry we refer to \cite{[P]} and \cite{[BBI], BGP} respectively.

\textbf{Acknowledgements.} The author would like to thank professor Karsten Grove for his enlightening suggestions and discussions. His suggestions on Gromov-Hausdorff precompact and Anderson-Cheeger compactness lead to Theorem B and C.

\section{Examples, Models and Metric Invariants}        

In this section we will cast the notion of \emph{tight cover} in terms of a metric invariant, exhibit examples,  and present model spaces, which in turn proves Theorem A.


The following notion captures how tight a compact (geodesic) metric space $(X, \rho)$ can be covered by two open balls in the sense described in the introduction.

\smallskip
Define the 2-\emph{covering defect} at $p,q \in X$ by
$$\cd(p,q) := \min \{d \ | \  D(p, r+d) \cup D(q, r+d) = X,  \  \rho(p,q) = 2r \}$$
and the  \emph{covering defect} of $X$ as

$$\cd(X): = \min _{p,q}\cd(p,q)$$

Clearly, then for any $d > \cd(X)$ there is a pair of points $p,q \in X$ such that $X$ is covered by the open $d+ \rho(p,q)/2$-balls centered at $p$ and $q$. The covering assumption in the Theorems B, C and D is that its covering defect $\cd(M)$ is at most $\epsilon$. Note also, that the covering defect of $X$ satisfies $0 \le \cd(X) \le \diam X/2$, and $\cd(X) = 0$ if and only if $X$ is \emph{efficiently covered by two closed metric balls} in the following sense:
$$X = D(p,r) \cup D(q,r), \ \text{with} \  \rho(p,q) = 2r.$$
Moreover, the maximal value $\cd(X) = \diam X/2$ is achieved, e.g. for $X$ a projective space with its rank one symmetric metric.

\medskip

We note that by continuity, any Gromov-Hausdorff limit $X$ of a sequence of Riemannian manifolds $M_i$ with $\cd(M_i)$ approaching $0$ has $\cd(X) = 0$. This suggest many types of examples.


\begin{example}[Sphere with tiny surgery]\label{surgery}The unit $n$-sphere $\Sph^n$ has $\cd(\Sph^n) = 0$. Now any closed manifold $M$ has a Riemannian metric arbitrarily Gromov-Hausdorff close to $\Sph^n$: Simply take any Riemannian metric on $M - D$ which is a product near the boundary sphere (with constant curvature). Scale it to any small size, sew it to the complement of a suitable small disc in $\Sph^n$ and smooth out.
\end{example}


As a very different type of $X$ with $\cd(X) = 0$, consider, e.g., the wedge $X = \Disc^n \vee \Disc^n$ of the euclidean unit disc $\Disc^n$ with itself at a point of the boundary $\Sph^{n-1}$. Concretely we may take  $X \subset \R^n$ with centers of the discs say $(-1, 0,\ldots,0)$ and $(1,0,\ldots,0)$. Note that $\cd(X)$ is realized at all the pairs of points $(-s, 0,\ldots,0)$ and $(s,0,\ldots,0)$, $1 \le s \le 2$.

\begin{example} [Connected Sum]\label{sum}
Let $M$ and $N$ be closed Riemannian manifolds with diameter $r$, and $X = M \vee N$, where $M$ and $N$ are joined at a point where the diameter is realized. Clearly, $\cd(M \vee N) = 0$ and hence the connected sum $M\# N$ has a Riemannian metric with $\cd(M\# N)$ arbitrarily close to $0$.
\end{example}

A feature of the latter examples is that \emph{geodesics branch}, at least at the point where the two spaces are joined. We will say that a geodesic metric space $X$ is \emph{non-branching} if any minimal geodesics, say  $c_{xy}: [0,2\epsilon] \to X$,  $c_{xz}: [0,2\epsilon] \to X$ from $x$ to $y$, and $x$ to $z$ respectively, coincide if $c_{xy}(\epsilon) = c_{xz}(\epsilon)$.

\smallskip
Our results pivot around the following observation which combined with continuity of $\cd$ and $\text{excess}$ proves Theorem A in the introduction.

\begin{lem}[Suspension Rigidity]\label{suspension}
Let $X = D(p,r) \cup D(q,r)$ be a compact non-branching geodesic metric space with $\rho(p,q) = 2r$. Then

$\bullet$ $\partial D(p,r) = \partial D(q,r) : = E$, and

$\bullet$ $X$ is topologically the suspension of $E$

$\bullet$ $\excess X = 0$

\no In fact, there is a unique minimal geodesic from $p$ to $q$ through every $x\in X - \{p,q\}$, and every metric sphere $S(p,s) = \partial D(p,s) = \partial D(q,2r-s) = S(q,2r-s)$ with center $p$ (and $q$) is homeomorphic to $E$.
\end{lem}

\begin{proof}
We first show that $\partial D(p,r) = E = \partial D(q,r)$. Let $x \in \partial D(p,r)$, i.e., $\rho(p,x) = r$ and pick normal minimal geodesics $c_{px}$, respectively $c_{qx}$ from $p$, respectively $q$ to $x$. We need to show that $\rho(x,q) = r$. Since clearly  $\rho(x,q) \ge r$, suppose $\rho(x,q) > r$ and let $x' = c_{qx}(r)$. By assumption,  the restriction of $c_{qx}$ from $x'$ to $x$ is contained in $D(p,r)$, and it follows that $\rho(x',p) = r$. Now let $c_{px'}$ be a minimal geodesic from $p$ to $x'$. It follows that the concatenation of $c_{px'}$ with the restriction of $c_{qx}$ (in the opposite direction) from $x'$ to $q$ is a minimal geodesic from $p$ to $q$. Since geodesics do not branch, this however implies that $x \in B(p,r)$, a contradiction.

Again, since geodesics do not branch in $X$, we see that there is a unique minimal geodesic from any $x \in E$ to $p$ as well as to $q$. Now let $y \in B(p,r)$ and consider minimal geodesics $c_{py}$ and $c_{yq}$. Since $x= c_{qy}(r) \in E$, we see as above that  $c_{py}$ together with $c_{yq}$ form a minimal geodesic from $p$ to $q$.

All in all we have seen that any point $y \in X -\{p,q\}$ lies on a unique minimal geodesic of length $2r$ from $p$ to $q$, and hence $X$ is the suspension of $E \cong S(p,s)$ for any $0<s<2r$.
\end{proof}

\begin{rem}
Note, that the assumption $X = D(p,r) \cup D(q,r)$ can be replaced by $X = D(p,s) \cup D(q,d-s)$ for any fixed $s$ where $0<s<d = \rho(p,q)$, with the same conclusion.
\end{rem}

\begin{rem}[Pointed Wiedersehen's Manifold]
If $X=M$ is a Riemannian manifold satisfying the assumptions above, it follows that the cut locus of $p$ is $q$ (and vice versa). Clearly, such a manifold is homeomorphic to a sphere. Conversely, A. Weinstein showed that any twisted (exotic) sphere has a Riemannian metric of this kind (see \cite{Be},Appendix C.4).
\end{rem}


If $V_{p}(r)$ denote the volume of $B_{p}(r) \subset M$, we also have the immediate

\begin{cor}
A Riemannian manifold $M$ with $V(M)=V_{p}(s)+V_{q}(d-s)$ for some $0 < s< \rho(p,q) = d$ is a pointed wiedersehen's manifold.
\end{cor}

\begin{proof}
The volume assumption implies $M^{n}=\overline{B}_{p}(s)\cup \overline{B}_{q}(d-s)$.
\end{proof}

\medskip

Note that Theorems B and C are also direct consequence of lemma \ref{suspension} via a Gromov-Hausdorff limit argument once it is proved that the limit $X$ is indeed topologically an $n$-sphere. For Theorem B this follows from the fact that small metric spheres of the limit are indeed $n-1$ spheres (a consequence of \cite{AC} as proved in \cite{[Pet]}). For Theorem C this is a consequence of the observation that the space of directions $S_pX$ at the suspension points $p \in X$, consists of geodesic directions, and that this space is homeomorphic to the spheres centered at $p$. Moreover, from \cite{Ka2} we know that $S_pX$ is itself the Gromov-Hausdorff limits of $n-1$ spheres with a lower curvature bound, and hence by Perelman's Stability Theorem \cite{Ka1} $S_pX$ is homeomorphic to $\Sph^{n-1}$.


The following examples show that for fixed $n$ , the Theorems B and C are optimal in the sense that the conclusion fails if any of the assumptions are missing.

\begin{example} [Curvature needed]The sphere with tiny surgery \ref{surgery} provides examples with lower bound on volume and upper bound on diameter, but where the lower curvature bound necessarily goes to minus infinity.
\end{example}

\begin{example}[Volume needed] Let $M$ be a closed manifold supporting an isometric cohomogeneity one action by $\G$. Then $M/\G$ is either an interval or a circle, but in either case it satisfies the exact ball cover condition. Cheeger deforming $M$ results in metrics keeping a lower curvature bound \cite{Ch}, the diameter convergies to $\diam M/\G = d$, but the volume approaches zero. Clearly, for any $\epsilon$ a sufficiently Cheeger deformed $M$ will be covered by to balls of radius $d/2 + \epsilon$.
\end{example}

\begin{example}[Diameter needed] Consider a sphere $N$ with nonnegative curvature like a cylinder with two spherical caps, say with radius $r$. Let $M$ be the manifold obtained from $N$ by attaching a tine handle to one of the cap's. As the length $d$ of cylinder goes to infinity, $M$ is covered by balls of radius $d/2 + \epsilon$ with $\epsilon$ going to zero. Throughout, the lowner curvature bound stays the same, the volume goes to infinity, and the topology of $M$ stays constant.
\end{example}

Without any curvature assumption, recall that the class of closed $n$ manifolds with a lower injectivity radius bound and an upper volume bound is Gromov-Hausdorff pre compact (cf. \cite{GPW,Cr}), it is natural to pose the following

\begin{problem}
Given $n$, $i$ and $V$. Is there an $\epsilon_3 = \epsilon(n,i,V)$ such that any closed $n$-manifold $M$ with $\inj M \ge i$, $\vol M \le V$ and $\cd(M) < \epsilon_3$ is homeomorphic to $\Sph^n$ when $\epsilon < \epsilon_3$.
\end{problem}

A limiting object $X$ will have excess 0, but geodesics (if they exist) may branch and not much is known about the geometry of limit objects from this class. Topologically they are homology manifolds.

\section{A constructive approach}   

In this section we show that if the covering defect is attained at a pair of points at maximal distance in $M$, then an explicit estimate for $\epsilon_1$ in Theorem C can be given.

In contrast to the proofs in the previous section, we establish directly that the open metric balls $B(p,r+\epsilon)$ and $B(q,r+\epsilon)$ are homeomorphic to a euclidean ball and hence $M$ is a sphere by Brown's Theorem. This in turn is done using Toponogov's distance comparison Theorem to show that the distance functions $\rho_p$ and $\rho_q$ have no critical points in their respective balls in the sense of \cite{[GS]} (see also, e.g., the survey \cite{Gr}). In particular, under this scenario, the ``$\epsilon_1$" will not depend on $n=\dim M$.

\begin{thm}\label{explicit}
There is an explicit $\epsilon_2 = \epsilon_2(k,i,d)$, such that any closed Riemannian $n$-manifold $M$ with $\sec M \ge k$, $\inj M \ge i$, $\rho(p,q)= d = \diam M$ and $\cd(p,q) < \epsilon_2$ is homeomorphic to $n$-sphere $\Sph^{n}$.

For, e.g., $\sec M \ge 0$ the estimate $\epsilon_2(i,d) = \frac{1}{2}(-d+\sqrt{d^{2}+\frac{i^{2}}{2}})$ can be used.
\end{thm}

For simplicity, we provide an explicit estimate for $\epsilon_2$ when $\sec M \ge 0$, and outline an approach via an excess estimate (and \cite{[GP]}) how to derive it in the general case of $\sec M \ge k$ (see \ref{excess}).
\smallskip

\begin{proof}
Since $\inj M \ge i$, it suffices to prove that $\rho_p$ has no critical points in $B(p,d/2+ \epsilon) - B(p,i)$, the argument for $\rho_q$ being exactly the same.

\smallskip

For any $x \in B(p,d/2+ \epsilon) - B(p,i)$, let $c_{xq}: [0, l_2+l_4] \to M$ be a normal minimal geodesic from $x$ to $q$ with $y = c_{xq}(l_2) \in \partial B(p,d/2+ \epsilon)$ and $\rho(y,q) = l_4$. Moreover, let $c_{xp}: [0,l_1] \to M$  be a normal minimal geodesic making an angle $\theta$ with $c_{xq}$. Finally, let $c_{yp}: [0,l_3] \to M$ be a normal minimal geodesic from $y$ to $p$ making and angle $\alpha$ with $c_{xq}$ at $y$.

The overall strategy will be to derive an upper bound for $l_2$ using the assumption $\rho(x,q) = l_2+l_4 \le d$, and a lower bound for $l_2$ under the assumption that $x$ is a critical point for $\rho_p$. Moreover, these bounds will force a lower for $\cd(pq)$.

First we have the following initial length estimates:
\smallskip

\begin{enumerate}
\item
$i \le l_1 < \frac{d}{2}+\epsilon$.
\item
 $l_{2}<\frac{d}{2}+\epsilon$.
 \item
 $l_{3}=\frac{d}{2}+\epsilon$.
 \item
 $\frac{d}{2}-\epsilon<l_{4}<\frac{d}{2}+\epsilon$.
 \end{enumerate}
Here (1) and (3) are immediate. By the covering condition, we know that $y\in B_{q}(\frac{d}{2}+\epsilon)$, and hence the upper bound in (4). Also by assumption $l_{2}+l_{4}=\rho(x,q)\leq d$, and hence $l_{2}\leq d-l_{4}$. Now from the triangle $l_{4} \ge d-l_{3}=\frac{d}{2}-\epsilon$ with equality only if $x$ is an a minimal geodesic from $p$ to $q$. The latter, however is clearly impossible, establishing the left inequality of (4) and in turn (2).

\medskip

The hinge version of the Toponogov's comparison theorem applied to $\Delta pyq$ yields:
$$d^{2} \leq  l_{3}^{2}+l_{4}^{2}-2l_{3}l_{4}\cos\alpha, \ \ \text{and hence} \  \  \cos\alpha \leq\frac{l_{3}^{2}+l_{4}^{2}-d^{2}}{2l_{3}l_{4}}.$$

Since $d>l_{3}$, $\frac{l_{3}^{2}+l_{4}^{2}-d^{2}}{2l_{3}l_{4}}$ is an increasing function of $l_{4}$, so
$$\frac{l_{3}^{2}+l_{4}^{2}-d^{2}}{2l_{3}l_{4}}<\frac{(\frac{d}{2}+\epsilon)^{2}+
(\frac{d}{2}+\epsilon)^{2}-d^{2}}{2(\frac{d}{2}+\epsilon)^{2}}=1-\frac{d^{2}}{2(\frac{d}{2}+\epsilon)^{2}}.$$

This then leads to the following estimate
\begin{equation}\label{alpha}
\cos\alpha<1-\frac{d^{2}}{2(\frac{d}{2}+\epsilon)^{2}}.
\end{equation}

Applying the hinge version of Toponogov's  comparison theorem to $\Delta pyx$ we have:

$$l_{1}^{2}\leq l_{2}^{2}+l_{3}^{2}-2l_{2}l_{3}\cos(\pi-\alpha)=l_{2}^{2}+l_{3}^{2}+2l_{2}l_{3}\cos\alpha.$$

\no or equivalently
$$l_{1}^{2}\leq (l_{3}-l_{2})^{2}+2l_{2}l_{3}(\cos\alpha+1).$$

\smallskip
\no But from (1), (2) and the estimate \ref{alpha} on $\alpha$ above, we estimate that last term by

$$l_{2}l_{3}(\cos\alpha+1)< (\frac{d}{2}+\epsilon)^{2}[2-\frac{d^{2}}{2(\frac{d}{2}+\epsilon)^{2}}] = 2\epsilon(d+\epsilon).$$

\no Thus by (1)
$$i^{2}\leq l_{1}^{2}< (l_{3}-l_{2})^{2}+4\epsilon(d+\epsilon).$$

\no and hence our desired upper estimate for $l_2$:
\begin{equation}\label{upper}
l_{2} < l_{3}-\sqrt{i^{2}-4\epsilon(d+\epsilon)} =
 \frac{d}{2}+\epsilon-\sqrt{i^{2}-4\epsilon(d+\epsilon)}.
\end{equation}

\medskip

We now proceed to derive the desired lower estimate for $l_2$ under the assumption that $\theta \le \pi/2$, which can be achieved if $x$ is a critical point for $\rho_p$.

In this case, the hinge version of Toponogov's comparison theorem yields

$$l_{3}^{2}\leq l_{1}^{2}+l_{2}^{2}.$$

\no Combining this with the hinge estimate on ${l_1}^2$ above we get

$$0\leq 2l_{2}^{2}+2l_{2}l_{3}\cos\alpha.$$

\no Substituting the estimate \ref{alpha} on $\alpha$ above, we get

$$0< l_{2}+(\frac{d}{2}+\epsilon)[1-\frac{d^{2}}{2(\frac{d}{2}+\epsilon)^{2}}],$$

\no leading to the desired lower estimate on $l_2$:
\begin{equation}\label{lower}
(\frac{d}{2}+\epsilon)[\frac{d^{2}}{2(\frac{d}{2}+\epsilon)^{2}}-1]  <  l_{2}.
\end{equation}

\smallskip

Clearly, the lower and upper estimates for $l_2$ contradict one another when $\epsilon$ is small enough. In this case, therefore $x$ cannot be a critical point for $\rho_p$.

We proceed to find an explicit estimate for this by combining \ref{lower} and \ref{upper}:

$$(\frac{d}{2}+\epsilon)[\frac{d^{2}}{2(\frac{d}{2}+\epsilon)^{2}}-1]\geq\frac{d}{2}+\epsilon-\sqrt{i^{2}-4\epsilon(d+\epsilon)}.$$

This is equivalent to

\begin{equation}
8(d\epsilon+\epsilon^{2})-i^{2}+\frac{d^{4}}{4(\frac{d}{2}+\epsilon)^{2}}-d^{2}\leq0.
\end{equation}

Because we always have $\frac{d^{4}}{4(\frac{d}{2}+\epsilon)^{2}}-d^{2}<0$, it suffices to make sure that $8(d\epsilon+\epsilon^{2})-i^{2}\leq0$. Thus  we can choose
$$\epsilon(d, i)=\frac{1}{2}(-d+\sqrt{d^{2}+i^{2}/2}).$$
Indeed, when $\varepsilon<\epsilon(d, i)$, $\rho_p$ cannot have any critical points even in $D(p,d/2 + \epsilon)$.
\end{proof}

As indicated, to complete the proof Theorem 2.1 for general $k \in R$, we proceed via Lemma
3 in \cite{[GP]} and the following

\begin{thm}\label{excess}
Any closed Riemannian $n$-manifold $M$ with $\sec M \ge k$, $\rho(p,q)= d = \diam M$ and $\cd(p,q) < \epsilon$ has $\excess M < f(\epsilon, k, d)$ with $f(\epsilon,k,d)$ going to $0$ as $\epsilon$ goes to $0$.
\end{thm}

\begin{proof}
By scaling, it is enough to consider the case $\sec M \ge -1$. Also, it suffices to estimate the excess $\rho(p,x) + \rho(x,q) - \rho(p,q)$  for any $x \in B(p,d/2 + \epsilon)$. Using the notation from the proof of Theorem \ref{explicit}, this is
$$l_1+l_2+l_4 - d = (l_1+l_2 - l_3) + (l_3 + l_4 - d),$$
where the key to estimate each summand is an estimate on the angle $\alpha$. Clearly, $\alpha$ exceeds the angle, $\alpha_0$ in constant curvature $-1$ of a triangle with two sides of length $l_3 = d/2 +\epsilon$ and base of length $d = \rho(p,q)$, i.e.,
$$\cosh d = \cosh^2(d/2 + \epsilon) -\sinh^2(d/2+ \epsilon)\cos \alpha_0.$$
Also, $l_1+l_2 - l_3$ is maximized in the limit when $l_2 = d/2 + \epsilon$. Thus
$$\rho(p,x) + \rho(x,q) - \rho(p,q) \le 2 \epsilon + l_0,$$
where
$$\cosh l_0 = \cosh^2(d/2 + \epsilon) - \sinh^2(d/2 + \epsilon)\cos (\pi-\alpha_0)$$.

Since $\alpha_0$ approaches $\pi$ for any fixed $d$ as $\epsilon$ approaches $0$, and likewise $l_0$ approaches $0$ for any fixed $d$ and $\epsilon$ as $\pi - \alpha_0$ approaches $0$, this completes the proof.
\end{proof}

\medskip

In the special case of even dimensional, simply connected manifolds $M$ with positive curvature, say scaled as $0< \sec M \le 1$, recall that the injectivity radius is bounded below as $\inj M \ge \pi$ \cite{Kl}. Thus

\begin{cor}
Any orientable closed Riemannian manifold $M^{2n}$ with $0<\sec_{M}\leq 1$ and
$$M^{n}=B_{p}(\frac{d}{2}+\epsilon)\cup B_{q}(\frac{d}{2}+\epsilon), \ where \   \rho(p,q) = d = \diam M$$
is homeomorphic to $\Sph^{2n}$ as long as   $\epsilon < \frac{1}{2}(\sqrt{d^{2}+\pi^{2}/2} -d).$

\end{cor}


\providecommand{\bysame}{\leavevmode\hbox    
to3em{\hrulefill}\thinspace}                                            

\end{document}